\documentclass{amsart}

\usepackage{amsmath,amssymb,amsthm}
\usepackage[all]{xy}

\DeclareMathOperator{\Diff}{Diff}
\DeclareMathOperator{\Aut}{Aut}
\DeclareMathOperator{\Lie}{Lie}

\DeclareMathOperator{\SL}{SL}

\DeclareMathOperator{\GCD}{GCD}
\DeclareMathOperator{\hol}{hol}
\DeclareMathOperator{\id}{id}
\DeclareMathOperator{\vol}{vol}
\DeclareMathOperator{\dev}{dev}

\DeclareMathOperator{\modular}{mod}

\newtheorem{thm}{Theorem}
\newtheorem{lem}[thm]{Lemma}
\newtheorem{prop}[thm]{Proposition}
\newtheorem{cor}[thm]{Corollary}
\newtheorem{defn}[thm]{Definition}

\begin{document}
\title{Rigidity of the \'{A}lvarez class}
\author{Hiraku Nozawa}
\address{Unit\'{e} de Math\'{e}matiques Pures et Appliqu\'{e}es, \'{E}cole Normale Sup\'{e}rieure de Lyon, \newline
46, all\'{e}e d'Italie 69364 Lyon Cedex 07, France.}
\email{nozawahiraku@06.alumni.u-tokyo.ac.jp}

\begin{abstract}
Let $(M,\mathcal{F})$ be a closed manifold with a Riemannian foliation. The \'{A}lvarez class of $(M,\mathcal{F})$ is a cohomology class of $M$ of degree $1$ whose triviality characterizes the minimizability or the geometrically tautness of $(M,\mathcal{F})$. We show that the integral of the \'{A}lvarez class of $(M,\mathcal{F})$ along every closed path is the logarithm of an algebraic integer if $\pi_{1}M$ is polycyclic or $\mathcal{F}$ is of polynomial growth.
\end{abstract}


\maketitle

\section{Introduction}

A foliated manifold $(M,\mathcal{F})$ is minimizable or geometrically taut if there exists a Riemannian metric $g$ on $M$ such that every leaf of $\mathcal{F}$ is a minimal submanifold of $(M,g)$. \'{A}lvarez L\'{o}pez \cite{Alvarez Lopez} defined a cohomology class $\xi(\mathcal{F})$ of degree $1$ for a closed manifold $M$ with a Riemannian foliation $\mathcal{F}$ whose triviality is equivalent to the minimizability of $(M,\mathcal{F})$. Here, $\xi(\mathcal{F})$ is called the \'{A}lvarez class of $(M,\mathcal{F})$. \'{A}lvarez L\'{o}pez \cite{Alvarez Lopez 2} showed that $\xi(\mathcal{F})$ is equal to a secondary characteristic class of the Molino's commuting sheaf of $(M,\mathcal{F})$ up to multiplication of a non-zero constant. These results of \'{A}lvarez L\'{o}pez describe topological aspects of the minimizability of Riemannian foliations as well as the results of Ghys \cite{Ghys} and Masa \cite{Masa}. It is remarkable because the minimizability of general foliations is characterized in terms of dynamical tools (see Sullivan \cite{Sullivan} and Haefliger \cite{Haefliger}), which are not directly related with the topology of manifolds. In this paper, we show that every period of the \'{A}lvarez class is the logarithm of an algebraic integer under certain topological conditions on $(M,\mathcal{F})$. Since the \'{A}lvarez class changes continuously for a smooth family of Riemannian foliations on a closed manifold, as shown in \cite{Nozawa}, the totally disconnectedness of the period of the \'{A}lvarez classes shown in this paper implies that they are invariant under deformation. This rigidity phenomenon of the \'{A}lvarez class is different from the known rigidity of secondary characteristic classes of flat vector bundles (see \cite{Kamber Tondeur}). 

We state the main result in this paper. We refer to \cite{Molino} or \cite{Moerdijk Mrcun} for the Riemannian and Lie foliation theories used in this paper. By the Molino's structure theorem, we mean Theorem 4.26 of \cite{Moerdijk Mrcun}. Let $M$ be a closed manifold and $\mathcal{F}$ a Riemannian foliation on $M$. Let $M^1$ be the transverse orthonormal frame bundle of $(M,\mathcal{F})$ for a bundle-like metric on $(M,\mathcal{F})$. Let $\mathcal{F}^{1}$ be the horizontal lift of $\mathcal{F}$. Let $F$ be the closure of a leaf of $(M^{1},\mathcal{F}^{1})$. Then, by the Molino's structure theorem, $(F,\mathcal{F}^{1}|_{F})$ is a $G$-Lie foliation for some Lie group $G$. Fix a point $x_{0}$ on $F$. Let $\hol_{F} \colon \pi_{1}(F,x_{0}) \longrightarrow G$ be the holonomy homomorphism of the Lie foliation $\mathcal{F}^{1}|_{F}$. The holonomy group $\hol_{F}(\pi_{1}(F,x_{0}))$ is independent of the bundle-like metric. Our main result in this paper is the following: 

\begin{thm}\label{Theorem : Rigidity}
If $\hol_{F}(\pi_{1}(F,x_{0}))$ is polycyclic, then the integral of the \'{A}lvarez class of $(M,\mathcal{F})$ along every closed path is the logarithm of an algebraic integer.
\end{thm}
\noindent See Definition \ref{Definition : Polycyclic groups} for the definition of polycyclic groups. Our main step of the proof of Theorem \ref{Theorem : Rigidity} is to establish a relation between the holonomy homomorphism of the Molino's commuting sheaf with $\hol_{F}(\pi_{1}(F,x_{0}))$ proved in Section 2. These fundamental objects in Riemannian foliation theory are related in this paper for the first time based on Lie's theorem. After the main step, we complete the proof by a rigid property of polycyclic groups and \'{A}lvarez L\'{o}pez's interpretation \cite{Alvarez Lopez 2} of the \'{A}lvarez class as a secondary characteristic class of Molino's commuting sheaf. Note that Meigniez \cite{Meigniez} constructed examples of Riemannian foliations on closed manifolds with non polycyclic fundamental groups whose \'{A}lvarez classes do not have the rigidity stated in Theorem \ref{Theorem : Rigidity} (see also \cite{Meigniez 2}).

As a direct result of the relation of the holonomy homomorphism of the Molino's commuting sheaf to $\hol_{F}(\pi_{1}(F,x))$ and $G$, we obtained a sufficient condition for the triviality of the \'{A}lvarez class. Then we get the following from a characterization of \'{A}lvarez L\'{o}pez \cite{Alvarez Lopez}:
\begin{thm}\label{Theorem : Reductive}
A Riemannian foliation on a closed manifold with semisimple structural Lie algebra is minimizable.
\end{thm}
\noindent The structural Lie algebra of $(M,\mathcal{F})$ is the Lie algebra of $G$ in the above notation.

We state corollaries of Theorem \ref{Theorem : Rigidity}. Let $(M,\mathcal{F})$ be a closed manifold with a Riemannian foliation. 
\begin{cor}\label{Corollary : Polycyclic}
If $\pi_{1}M$ is polycyclic, then the integral of the \'{A}lvarez class of $(M,\mathcal{F})$ along every closed path is the logarithm of an algebraic integer.
\end{cor}
\begin{cor}\label{Corollary : Nilpotent}
If $\mathcal{F}$ is of polynomial growth, then the integral of the \'{A}lvarez class of $(M,\mathcal{F})$ along every closed path is the logarithm of an algebraic integer.
\end{cor}
\noindent A foliated manifold $(M,\mathcal{F})$ is called of polynomial growth if every leaf of $\mathcal{F}$ is of polynomial growth. To deduce Corollary \ref{Corollary : Nilpotent} from Theorem \ref{Theorem : Rigidity}, we use a theorem of Carri\`{e}re \cite{Carriere} which claims that a Riemannian foliation is of polynomial growth if and only if the structural Lie algebra is nilpotent. A $1$-dimensional Riemannian foliation is called a Riemannian flow. Since every Riemannian flow is of polynomial growth, we have

\begin{cor}\label{Corollary : Riemannian Flow}
The integral of the \'{A}lvarez class of a Riemannian flow on a closed manifold along every closed path is the logarithm of an algebraic integer.
\end{cor}

The results in this paper were obtained during the author's stay in the University of Santiago de Compostela. The author is deeply grateful to Jes\'{u}s Antonio \'{A}lvarez L\'{o}pez for his invitation, great hospitality and valuable discussion. The author would like to express his gratitude to Yoshifumi Matsuda and Hiroki Kodama for pointing mistakes in the previous version of this paper. The author would like to express his gratitude to Shuji Yamamoto and Takeshi Katsura for the simplification of the proof of Lemma \ref{Lemma : Polycyclic groups 2}.

\section{Holonomy of Molino's commuting sheaves of transversely parallelizable foliations}

Let $(M,\mathcal{F})$ be a closed manifold with a transversely parallelizable foliation. We will relate the holonomy homomorphism of $(M,\mathcal{F})$ to the structure Lie group of $(M,\mathcal{F})$ and the holonomy group of the Lie foliation defined on each leaf closure. Lie's theorem is the key. We prepare our notation. By the Molino's structure theorem \cite{Molino}, $M$ is the total space of a smooth fiber bundle $p \colon M \longrightarrow W$ such that each fiber of $p$ is the closure of a leaf of $\mathcal{F}$. We fix a fiber $F$ of $p$. Then $p$ is a smooth $(F,\mathcal{F}|_{F})$-bundle by the Molino's structure theorem \cite{Molino}. 

We fix a point $x_{0}$ in $F$. By the Molino's structure theorem \cite{Molino}, $(F,\mathcal{F}|_{F})$ is a $G$-Lie foliation for some connected and simply-connected Lie group $G$. Let
\begin{equation}
\hol_{F} \colon \pi_{1}(F,x_{0}) \longrightarrow G
\end{equation}
be the holonomy homomorphism of the Lie foliation $\mathcal{F}|_{F}$. We put
\begin{equation}
\Gamma=\hol_{F}(\pi_{1}(F,x_{0}))
\end{equation}
and
\begin{equation}
\Aut(G,\Gamma)= \{ a \in \Aut(G) \ | \ a(\Gamma)=\Gamma \}.
\end{equation}

We refer to pages 125--130 of \cite{Molino} for the Molino's commuting sheaf $\mathcal{C}$ of $(M,\mathcal{F})$. We summarize its definition and properties. For an open set $U$ in $M$, $C_{b}^{\infty}(TU/T(\mathcal{F}|_{U}))$ is the space of transverse fields of $(U,\mathcal{F}|_{U})$ defined by
\begin{equation}
C_{b}^{\infty}(TU/T(\mathcal{F}|_{U})) = \{ Y \in C^{\infty}(TU/T(\mathcal{F}|_{U})) \ | \ [Z,Y]=0, \forall Z \in C^{\infty}(T(\mathcal{F}|_{U})) \}.
\end{equation}
We define a vector space $\mathcal{C}(U)$ by the space of transverse fields which commutes with global transverse fields as follows: 
\begin{equation}
\mathcal{C}(U) = \{ X \in C^{\infty}_{b}(TU/T(\mathcal{F}|_{U})) \ | \ [X,Y|_{U}]=0, \forall Y \in C_{b}^{\infty}(TM/T\mathcal{F}) \}.
\end{equation}
Note that Lie brackets of vector fields on $M$ induces Lie brackets $[\cdot,\cdot] \colon C^{\infty}(T\mathcal{F}) \times C^{\infty}(TM/T\mathcal{F}) \longrightarrow C^{\infty}(TM/T\mathcal{F})$ and $[\cdot,\cdot] \colon C_{b}^{\infty}(TM/T\mathcal{F}) \times C_{b}^{\infty}(TM/T\mathcal{F}) \longrightarrow C_{b}^{\infty}(TM/T\mathcal{F})$. Since $\mathcal{C}(U)$ is closed under the Lie bracket, $\mathcal{C}$ is a presheaf of Lie algebras. The sheafication of $\mathcal{C}$ is called the Molino's commuting sheaf of $(M,\mathcal{F})$.

By the homogeneity of transversely parallelizable foliations, for each two points $x$, $y \in M$, there exists a diffeomorphism of $M$ which preserves $\mathcal{F}$ and maps $x$ to $y$ (see Theorem 4.8 of \cite{Moerdijk Mrcun}). Since $\mathcal{C}$ is determined by $\mathcal{F}$, the stalks of $\mathcal{C}$ at $x$ and $y$ are isomorphic. Let $Y$ be an element of $C_{b}^{\infty}(TM/T\mathcal{F})$. The lift of $Y$ to $C^{\infty}(TM)$ generates an $\mathcal{F}$-preserving flow on $M$ (see Proposition 2.2 of \cite{Molino}). The orbits of the flow generated by $Y$ are well-defined in $M/\mathcal{F}$ (see 129 page of \cite{Molino}). If an element $s$ of $\mathcal{C}(U)$ satisfies $s(x)=0$ for a point $x$ in $U$, then $s$ is zero on the orbit of $x$ of the flow generated by $Y$ by the equation $[s,Y|_{U}]=0$. By the transversely parallelizability of $(M,\mathcal{F})$, the union of the orbits of $x$ of the flows generated by elements of $C_{b}^{\infty}(TM/T\mathcal{F})$ contains an open neighborhood of $x$. Hence if $s$ is zero at $x$, then $s$ is zero on an open neighborhood of $x$ in $M$. In other words, $\mathcal{C}$ is Hausdorff. Hence $\mathcal{C}$ is locally constant. We denote the holonomy homomorphism of the Molino's commuting sheaf as a locally constant sheaf on $M$ by 
\begin{equation}
m \colon \pi_{1}(M,x_{0}) \longrightarrow \Aut(\mathcal{C}_{x_{0}}).
\end{equation}

By Proposition 4.4 of \cite{Molino}, the stalk $\mathcal{C}_{x}$ of $\mathcal{C}$ at $x$ is identified with the Lie algebra of right invariant vector fields on $G$. This identification is given by the restriction map $\mathcal{C}_{x} \longrightarrow T_{x}M/T_{x}\mathcal{F}$ at $x$. Every element of $\mathcal{C}_{x_{0}}$ is tangent to the closures of the leaves of $\mathcal{F}$ and the image of the restriction map $\mathcal{C}_{x_{0}} \longrightarrow T_{x_{0}}M/T_{x_{0}}\mathcal{F}$ is $T_{x_{0}}F/T_{x_{0}}\mathcal{F}$ (see Theorem 4.3 of \cite{Molino}). Note that $T_{x_{0}}F/T_{x_{0}}\mathcal{F}$ is identified with the tangent space of $G$ at a point which is identified with the set of right invariant vector fields on $G$ by the definition of $G$-Lie foliations. 

We show
\begin{prop}\label{Proposition : Diagram}
We have the following commutative diagram of groups:
\begin{equation}\label{Equation ; Diagram}
\xymatrix{ \pi_{1}(M,x_{0}) \ar[r]^{m} \ar[d]_{h_{p}} & \Aut(\mathcal{C}_{x_{0}}) \\ 
\pi_{0}(\Diff(F,x_{0},\mathcal{F}|_{F})) \ar[ru]^{T} \ar[r]_>>>>>{\Phi} & \Aut(G,\Gamma) \ar[u]_{A}}
\end{equation}
where $\pi_{0}(\Diff(F,x_{0},\mathcal{F}|_{F}))$ is the group of the isotopy classes of the group of diffeomorphisms of $F$ which fix $x_{0}$ and map each leaf of $\mathcal{F}|_{F}$ to another leaf of $\mathcal{F}|_{F}$. Moreover $A$ is the action defined by $A(a)=Da$, where $Da$ is the restriction of the tangent map of $a \colon G \longrightarrow G$ to the Lie algebra of right invariant vector fields on $G$, which is identified to $\mathcal{C}_{x_{0}}$.
\end{prop}

\begin{proof}
We define the homomorphism $h_{p} \colon \pi_{1}(M,x_{0}) \longrightarrow \pi_{0}(\Diff(F,x_{0},\mathcal{F}|_{F}))$. Let $l$ be an element of $\pi_{1}(M,x_{0})$. We choose a representative $\gamma$ of $l$ so that $\gamma$ is a smooth map $[0,1] \longrightarrow M$ which induces a smooth map $\gamma_{S^1} \colon S^1 \longrightarrow M$, where we put $S^1=[0,1]/0 \sim 1$. Let $(\pi \circ \gamma_{S^1})^{*}\pi \colon (\pi \circ \gamma_{S^1})^{*}M \longrightarrow S^1$ be the pullback of the fiber bundle $\pi$ to $S^1$ by $\pi \circ \gamma_{S^1}$. Since $(\pi \circ \gamma_{S^1})^{*}M=\{(t,x) \in S^1 \times M \ | \ \pi \circ \gamma_{S^1}(t)=\pi(x) \}$ by the definition of the pullback of fiber bundles, we have a section 
\begin{equation}
\begin{array}{cccc}
s \colon & S^1 & \longrightarrow & (\pi \circ \gamma_{S^1})^{*}M \\
         &  t  & \longmapsto     & (t,\gamma(t)). 
\end{array}
\end{equation}
Considering the image of $s$ as a family of fiberwise base points, we regard $(\pi \circ \gamma_{S^1})^{*}\pi$ as an $(F,\mathcal{F}|_{F})$-bundle over $S^1$ with fiberwise base points. We can consider the pullback $(\pi \circ \gamma)^{*}\pi$ of $\pi$ to $[0,1]$ by $\pi \circ \gamma$ as an $(F,\mathcal{F}|_{F})$-bundle over $[0,1]$ with fiberwise base points in the same way. Since every Lie foliation is transversely parallelizable and every transversely parallelizable foliation is homogeneous (see Proposition 4.21 and Theorem 4.8 of \cite{Moerdijk Mrcun}), $(F,\mathcal{F}|_{F})$ is homogeneous; that is, the $\mathcal{F}$-preserving diffeomorphism group of $F$ acts on $F$ transitively. Hence we can trivialize $(\pi \circ \gamma)^{*}\pi$ as an $(F,\mathcal{F}|_{F})$-bundle with fiberwise base points. In other words, we have a trivialization $H \colon ((\pi \circ \gamma)^{*}M,(\pi \circ \gamma)^{*}\mathcal{F}) \longrightarrow ([0,1] \times F, [0,1] \times \mathcal{F}|_{F})$ which maps $\{ (t,\gamma(t)) \in (\pi \circ \gamma)^{*}M \ | \ t \in [0,1] \}$ to $\{(t,x_{0}) \in [0,1] \times F \ | \ t \in [0,1] \}$. Here $[0,1] \times \mathcal{F}|_{F}$ is a foliation on $[0,1] \times F$ with leaves $\{t\} \times L$ for $t \in [0,1]$ and $L \in \mathcal{F}|_{F}$. Hence there exists an element $f_{l}$ of $\Diff(F,x_{0},\mathcal{F}|_{F})$ such that $(\pi \circ \gamma_{S^1})^{*}\pi$ is isomorphic to the mapping torus $([0,1] \times F, [0,1] \times \mathcal{F}|_{F})/(0,f_{l}(x)) \sim (1,x)$ of $f_{l}$. We define $h_{p}(l)$ to be the isotopy class of $f_{l}$ in $\Diff(F,x_{0},\mathcal{F}|_{F})$.

We define $\Phi \colon \pi_{0}(\Diff(F,x_{0},\mathcal{F}|_{F})) \longrightarrow \Aut(G,\Gamma)$. Let $f \in \Diff(F,x_{0},\mathcal{F}|_{F})$. Let $X$ be a transverse field on $(F,\mathcal{F}|_{F})$. For any vector field $Y$ tangent to $\mathcal{F}|_{F}$ on $F$, we have $[Y,f_{*}X]=f_{*}[(f^{-1})_{*}Y,X]=0$. Then $f_{*}X$ is a transverse field. Hence $f_{*}$ induces an automorphism on the Lie algebra of transverse fields on $(F,\mathcal{F}|_{F})$. Since $\mathcal{F}|_{F}$ is dense in $F$, the Lie algebra of the transverse fields on $(F,\mathcal{F}|_{F})$ coincides with the structural Lie algebra $\Lie(G)$ of the $G$-Lie foliation $\mathcal{F}|_{F}$. Hence $f_{*}$ induces an automorphism of $\Lie(G)$. Hence, by Lie's theorem, there exists an automorphism $f_{G}$ of $G$ which induces $f_{*}$ on $\Lie(G)$. We define $\Phi(f)=f_{G}$. We show that $f_{G}$ maps $\hol_{F}(\pi_{1}(F,x_{0}))$ to itself. Let $\tilde{F}$ be the universal covering of $F$. We fix a point $\tilde{x}_{0}$ on the fiber of $x_{0}$ in $\tilde{F}$. Let $\dev \colon \tilde{F} \longrightarrow G$ be the developing map which maps $\tilde{x}_{0}$ to the unit element $e$ of $G$. Recall that we have $\dev \circ \tilde{\gamma} (1) = \hol_{F}([\gamma]) $ for an element $[\gamma]$ of $\pi_{1}(F,x_{0})$, where $\tilde{\gamma}$ is the lift of $\gamma$ to $\tilde{F}$ with $\tilde{\gamma}(0)=\tilde{x}_{0}$ (see Section 4.3.2 of \cite{Moerdijk Mrcun}). We denote the lift of $f$ on $\tilde{F}$ with $\tilde{f}(\tilde{x}_{0})=\tilde{x}_{0}$ by $\tilde{f}$. Let $\mathcal{G}$ be the pullback of $\mathcal{F}|_{F}$ to $\tilde{F}$. Clearly $\tilde{f}$ preserves $\mathcal{G}$. The leaves of $(\tilde{F},\mathcal{G})$ are the fibers of $\dev \colon \tilde{F} \longrightarrow G$ (see Lemma 4.23 of \cite{Moerdijk Mrcun}). Hence $\tilde{f}$ induces a diffeomorphism $\tilde{f}_{G}$ of $G$. Since $\tilde{f}_{G}$ and $f_{G}$ induces the same automorphism $f_{*}$ on $\Lie(G)$, we have $\tilde{f}_{G}=f_{G}$ by the uniqueness of Lie's theorem. Let $\widetilde{f \circ \gamma}$ be the lift of $f \circ \gamma$ to $\tilde{F}$ with $\widetilde{f \circ \gamma}(0)=\tilde{x}_{0}$. We have
\begin{equation}
\begin{array}{rl}
f_{G} \circ \hol_{F}([\gamma]) & = f_{G} \circ \dev \circ \tilde{\gamma} (1) \\
                               & = \tilde{f}_{G} \circ \dev \circ \tilde{\gamma} (1) \\
                               & = \dev \circ \tilde{f} \circ \tilde{\gamma} (1) \\
                               & = \dev \circ \widetilde{f \circ \gamma} (1) \\
                               & = \hol_{F}(f_{*} [\gamma]).
\end{array}
\end{equation}
Hence $f_{G}$ preserves the image of $\hol_{F}$. Since $\Gamma$ is finitely generated and dense in $G$, $\Aut(G,\Gamma)$ is totally disconnected in $\Aut(G)$. Then $\Phi(f)$ depends only on the isotopy class of $f$. 

We define $T \colon \pi_{0}(\Diff(F,x_{0},\mathcal{F}|_{F})) \longrightarrow \Aut(\mathcal{C}_{x_{0}})$. Let $f$ be an element of $\Diff(F,x_{0},\mathcal{F}|_{F})$. The tangent map $f_{*}$ of $f$ induces an automorphism on $T_{x_{0}}F/T_{x_{0}}\mathcal{F}$. $T_{x_{0}}F/T_{x_{0}}\mathcal{F}$ is identified with $\mathcal{C}_{x_{0}}$ in the way described in the paragraph previous to the statement of Proposition \ref{Proposition : Diagram}. We define $T(f)$ as an automorphism on $\mathcal{C}_{x_{0}}$ induced by $f_{*}$ through the identification with $T_{x_{0}}F/T_{x_{0}}\mathcal{F}$. By the last sentence of the previous paragraph, the automorphism on $T_{x_{0}}F/T_{x_{0}}\mathcal{F}$ induced by an element $f$ of $\Diff(F,x_{0},\mathcal{F}|_{F})$ depends only on the isotopy class of $f$. Hence $T \colon \pi_{0}(\Diff(F,x_{0},\mathcal{F}|_{F})) \longrightarrow \Aut(\mathcal{C}_{x_{0}})$ is defined.

We show the commutativity of the upper triangle of \eqref{Equation ; Diagram}. Since the Molino's commuting sheaf of $([0,1] \times F, [0,1] \times \mathcal{F}|_{F})$ is trivial, the map induced by $f$ on the stalk $\mathcal{C}_{x_{0}}$ of the Molino's commuting sheaf at $x_{0}$ is the holonomy of $\mathcal{C}$ with respect to $l$ by the definition of the holonomy homomorphism of locally constant sheaves. By the definition of $T$, it follows that the upper triangle commutes.

We show the commutativity of the lower triangle of \eqref{Equation ; Diagram}. For $f$ in $\Diff(F,x_{0},\mathcal{F}|_{F})$, both of $T(f)$ and $\Phi \circ A$ are maps induced on $T_{x_{0}}F/T_{x_{0}}\mathcal{F}$ from the tangent maps of $f$ at $x_{0}$ under the identification of $\mathcal{C}_{x_{0}}$ with $T_{x_{0}}F/T_{x_{0}}\mathcal{F}$. Hence they are equal.
\end{proof}

\section{Proof of Theorem \ref{Theorem : Reductive}}

Let $(M,\mathcal{F})$ be a closed manifold with a $p$-dimensional Riemannian foliation. Fix a bundle-like metric $g$ on $(M,\mathcal{F})$. Let $M^{1}$ be the transverse orthonormal frame bundle of $(M,\mathcal{F})$. Let $\mathcal{F}^{1}$ be the horizontal lift of $\mathcal{F}$ to $M^{1}$. We denote the canonical projection $M^{1} \longrightarrow M$ by $\pi$. We take a bundle-like metric $g^{1}$ on $(M^{1},\mathcal{F}^{1})$ so that $\pi$ is a Riemannian submersion and $T\mathcal{F}^{1}$ is orthogonal to the fibers of $\pi$.

\begin{lem}\label{Lemma : Pullback}
\begin{equation}\label{Equation : Kappa0}
\pi^{*} \xi(\mathcal{F}) = \xi(\mathcal{F}^{1}).
\end{equation}
\end{lem}
\begin{proof}
By the definition of $g$ and $g^{1}$, we have 
\begin{equation}\label{Equation : Kappa1}
\pi^{*}\chi=\chi^{1}
\end{equation}
where $\chi$ and $\chi^{1}$ are the characteristic forms of $(M,\mathcal{F},g)$ and $(M^{1},\mathcal{F}^{1},g^{1})$ respectively. Recall that the characteristic form of a Riemannian manifold $(N,h)$ with a $k$-dimensional foliation $\mathcal{G}$ is the $k$-form $\omega$ which vanishes when some vector is orthogonal to $\mathcal{G}$, and such that its restriction to each leaf of $\mathcal{G}$ is a volume form of norm one. We put $Q=(T\mathcal{F})^{\perp}$ and $Q^{1}=(T\mathcal{F}^{1})^{\perp}$. Let $d_{1,0} \colon C^{\infty}(\wedge^{p}T^{*}M) \longrightarrow C^{\infty}(Q^{*}) \otimes C^{\infty}(\wedge^{p}T^{*}\mathcal{F})$ be defined by the composition of the de Rham differential and orthogonal projection $C^{\infty}(\wedge^{p+1}T^{*}M) \longrightarrow C^{\infty}(Q^{*}) \otimes C^{\infty}(\wedge^{p}T^{*}\mathcal{F})$. Similarly, $d_{1,0}^{1} \colon C^{\infty}(\wedge^{p}T^{*}M^{1}) \longrightarrow C^{\infty}(Q^{1 *}) \otimes C^{\infty}(\wedge^{p}T^{*}\mathcal{F}^{1})$ is defined for $(M^{1},\mathcal{F}^{1})$. By the Rummler's formula (see \cite{Rummler} and \cite{Alvarez Lopez}), we have
\begin{equation}\label{Equation : Kappa2}
d_{1,0}\chi=-\kappa \wedge \chi, \quad d_{1,0}^{1}\chi^{1}=-\kappa^{1} \wedge \chi^{1}
\end{equation}
where $\kappa$ and $\kappa^{1}$ are the mean curvature forms of $(M,\mathcal{F},g)$ and $(M^{1},\mathcal{F}^{1},g^{1})$, respectively (see Chapter 10 of \cite{Candel Conlon}).  Note that $\kappa$ is an element of $C^{\infty}(Q^{*})$ and $\kappa^{1}$ is an element of $C^{\infty}(Q^{1 *})$ by the Rummler's formula. Since $\pi_{*}$ maps $Q^{1}$ to $Q$ and $T\mathcal{F}^{1}$ to $\mathcal{F}$, we have
\begin{equation}\label{Equation : Kappa3}
d_{1,0}^{1} \pi^{*}=\pi^{*} d_{1,0}
\end{equation}
By \eqref{Equation : Kappa1}, \eqref{Equation : Kappa2} and \eqref{Equation : Kappa3}, we have 
\begin{equation}\label{Equation : Kappa4}
\pi^{*} \kappa = \kappa^{1}.
\end{equation}
Let $p \colon M^{2} \longrightarrow M^{1}$ be the transverse orthonormal frame bundle of $(M^{1},\mathcal{F}^{1})$. Let $\mathcal{F}^{2}$ be the horizontal lift of $\mathcal{F}^{1}$ to $M^{2}$. We denote the foliations on $M^{i}$ which is defined by the closures of leaves of $\mathcal{F}^{i}$ by $\overline{\mathcal{F}}^{i}$. We denote the characteristic form of $\overline{\mathcal{F}}^{1}$ by $\overline{\chi}^{1}$. Since $p$ maps each leaf of $\mathcal{F}^{2}$ to each leaf of $\mathcal{F}^{1}$ diffeomorphically, we can define a characteristic form $\overline{\chi}^{2}$ of $\overline{\mathcal{F}}^{2}$ by
\begin{equation}\label{Equation : Kappa5}
\overline{\chi}^{2}=p^{*}\overline{\chi}^{1}.
\end{equation}
Let $\vol_{p}$ be a fiberwise volume form of $p$. Let $\vol_{\pi}$ be a fiberwise volume form of $\pi$. By the definition of the basic components of $\kappa$ and $\kappa^{1}$ (see \cite{Alvarez Lopez}), we have
\begin{equation}\label{Equation : Kappa6}
\kappa_{b} =  \frac{1}{\int_{\pi} \vol_{\pi}}  \int_{\pi} \rho_{\overline{\mathcal{F}}^{1}}( \pi^{*}\kappa ) \wedge \vol_{\pi} 
\end{equation}
and
\begin{equation}\label{Equation : Kappa7}
\kappa_{b}^{1} = \frac{1}{\int_{p} \vol_{p}} \int_{p} \rho_{\overline{\mathcal{F}}^{2}}( p^{*}\kappa^{1}) \wedge \vol_{p},
\end{equation}
where $\int_{\pi}$ and $\int_{p}$ are the integration along fibers and $\rho_{\overline{\mathcal{F}}^{i}}$ is the projection from the space of $1$-forms to the space of basic $1$-forms of $(M^{i},\mathcal{F}^{i})$ defined for $i=1$ and $2$ as follows:
\begin{equation}
\rho_{\overline{\mathcal{F}}^{i}} \Big(\sum_{j=1}^{q^{i}} f_{j} \omega_{j}^{i} + \alpha \Big) =  \frac{1}{\int_{\overline{\mathcal{F}}^{i}} \overline{\chi}^{i}} \sum_{j=1}^{q^{i}} \Big(\int_{\overline{\mathcal{F}}^{i}} f_{j} \overline{\chi}^{i} \Big) \omega_{j}^{i}
\end{equation}
for $f_{1}$, $f_{2}$, $\cdots$, $f_{q^{i}}$ in $C^{\infty}(M^{i})$ and $\alpha$ in $C^{\infty}(T^{*}\mathcal{F}^{i})$, where $q^{i}$ is the codimension of $\mathcal{F}^{i}$, $\{\omega_{j}^{i}\}_{j=1}^{q^{i}}$ is a transverse parallelism of $(M^{i},\mathcal{F}^{i})$ and $\int_{\overline{\mathcal{F}}^{i}}$ is the integration along fibers. By \eqref{Equation : Kappa5} and the fact that $p$ maps each leaf of $\overline{\mathcal{F}}^{2}$ to each leaf of $\overline{\mathcal{F}}^{1}$, we have
\begin{equation}
p^{*} \int_{\overline{\mathcal{F}}^{1}} f \overline{\chi}^{1} = \int_{\overline{\mathcal{F}}^{2}} (p^{*}f) \overline{\chi}^{2}
\end{equation}
for $f$ in $C^{\infty}(M^{1})$. Hence we have
\begin{equation}\label{Equation : Kappa8}
\rho_{\overline{\mathcal{F}}^{2}} p^{*} = p^{*} \rho_{\overline{\mathcal{F}}^{1}}.
\end{equation}
For every basic $1$-form $\alpha$ of $(M^{1},\mathcal{F}^{1})$ whose restriction to each fiber of $\pi$ vanishes, we have
\begin{equation}
\alpha = \pi^{*} \left( \frac{1}{\int_{\pi} \vol_{\pi}} \int_{\pi} \alpha \wedge \vol_{\pi} \right).
\end{equation}
Since $\rho_{\overline{\mathcal{F}}^{1}} \pi^{*} \kappa$ is basic and the restriction of $\rho_{\overline{\mathcal{F}}^{1}} \pi^{*} \kappa$ to the fibers of $\pi$ vanishes, we have 
\begin{equation}\label{Equation : Kappa9}
\rho_{\overline{\mathcal{F}}^{1}} \pi^{*} \kappa = \pi^{*} \left( \frac{1}{\int_{\pi} \vol_{\pi}} \int_{\pi} \rho_{\overline{\mathcal{F}}^{1}} \pi^{*} \kappa \wedge \vol_{\pi} \right).
\end{equation}
Note that
\begin{equation}\label{Equation : Kappa10}
\frac{1}{\int_{p} \vol_{p}} \int_{p} p^{*} \beta \wedge \vol_{p}=\beta
\end{equation}
for every $\beta$ in $C^{\infty}(T^{*}M^{1})$. By \eqref{Equation : Kappa7}, \eqref{Equation : Kappa4}, \eqref{Equation : Kappa8}, \eqref{Equation : Kappa9}, \eqref{Equation : Kappa6} and \eqref{Equation : Kappa10}, we have
\begin{equation}\label{Equation : Kappa11}
\begin{array}{rl}
\kappa_{b}^{1} & = \frac{1}{\int_{p} \vol_{p}} \int_{p} \rho_{\overline{\mathcal{F}}^{2}}( p^{*}\kappa^{1}) \wedge \vol_{p} \\
 & = \frac{1}{\int_{p} \vol_{p}} \int_{p} \rho_{\overline{\mathcal{F}}^{2}}( p^{*} \pi^{*} \kappa) \wedge \vol_{p} \\
 & = \frac{1}{\int_{p} \vol_{p}} \int_{p} (p^{*} \rho_{\overline{\mathcal{F}}^{1}} \pi^{*} \kappa) \wedge \vol_{p} \\
 & = \frac{1}{\int_{p} \vol_{p}} \int_{p} p^{*} \pi^{*} \left( \frac{1}{\int_{\pi} \vol_{\pi}} (\int_{\pi} \rho_{\overline{\mathcal{F}}^{1}} \pi^{*} \kappa \wedge \vol_{\pi}) \right) \wedge \vol_{p} \\
 & = \frac{1}{\int_{p} \vol_{p}} \int_{p} p^{*} \pi^{*} \kappa_{b} \wedge \vol_{p} \\
 & = \pi^{*} \kappa_{b}.
\end{array}
\end{equation}
By the definition of the \'{A}lvarez class (see \cite{Alvarez Lopez}), we have
\begin{equation}\label{Equation : Kappa12}
\xi(\mathcal{F}) = [\kappa_{b}], \quad \xi(\mathcal{F}^{1}) = [\kappa_{b}^{1}].
\end{equation}
By \eqref{Equation : Kappa11} and \eqref{Equation : Kappa12}, we have \eqref{Equation : Kappa0}.
\end{proof}

We state two theorems of \'{A}lvarez L\'{o}pez.
\begin{thm}[\'{A}lvarez L\'{o}pez \cite{Alvarez Lopez}]\label{Theorem : Alvarez Lopez}
For a closed manifold $M$ with a Riemannian foliation $\mathcal{F}$, the \'{A}lvarez class $\xi(\mathcal{F})$ of $(M,\mathcal{F})$ in $H^{1}(M;\mathbb{R})$ is well-defined. $\xi(\mathcal{F})$ is trivial if and only if $(M,\mathcal{F})$ is minimizable.
\end{thm}

\begin{thm}[\'{A}lvarez L\'{o}pez \cite{Alvarez Lopez 2}]\label{Theorem : Alvarez Lopez 2}
For a closed manifold $M$ with a Riemannian foliation $\mathcal{F}$, the holonomy homomorphism of the determinant bundle of the Molino's commuting sheaf of $(M,\mathcal{F})$ coincides with the homomorphism defined by
\begin{equation}
\begin{array}{ccc}
\pi_{1}M & \longrightarrow & \mathbb{R} \\
  l & \longmapsto & e^{\int_{l} \xi(\mathcal{F})}
\end{array}
\end{equation}
where $\xi(\mathcal{F})$ is the \'{A}lvarez class of $(M,\mathcal{F})$.
\end{thm}

We show Theorem \ref{Theorem : Reductive} by Proposition \ref{Proposition : Diagram} and the above two theorems of \'{A}lvarez L\'{o}pez.

\begin{proof}[End of the proof of Theorem \ref{Theorem : Reductive}]
Let $(M,\mathcal{F})$ be a closed manifold with a Riemannian foliation whose structural Lie algebra is semisimple. Let $M^{1}$ be the transverse orthonormal frame bundle of $(M,\mathcal{F})$ for a bundle-like metric on $(M,\mathcal{F})$. Let $\mathcal{F}^{1}$ be the lift of $\mathcal{F}$. Let $F$ be the closure of a leaf of $(M^{1},\mathcal{F}^{1})$. Then, by the Molino's structure theorem, $(F,\mathcal{F}^{1}|_{F})$ is a $G$-Lie foliation for some connected Lie group $G$. Since every automorphism of a connected semisimple Lie group is an inner automorphism and a semisimple Lie group has an biinvariant volume form (see Proposition 1.121 and Corollary 8.31 of \cite{Knapp}), the image of $A \colon \Aut(G) \longrightarrow \Aut(\Lie(G))$ is contained in $\SL(\Lie(G))$. By Proposition \ref{Proposition : Diagram}, the image of $m$ is contained in $\SL(\mathcal{C}_{x_{0}})$. By Theorem \ref{Theorem : Alvarez Lopez 2} of \'{A}lvarez L\'{o}pez, the \'{A}lvarez class of $(M^{1},\mathcal{F}^{1})$ vanishes if and only if the image of $m$ is contained in $\SL(\mathcal{C}_{x_{0}})$. Hence the \'{A}lvarez class of $(M^{1},\mathcal{F}^{1})$ vanishes. Since $\pi_{*} \colon \pi_{1}M^{1} \longrightarrow \pi_{1}M$ is surjective by the homotopy exact sequence of fiber bundles, $\pi^{*} \colon H^{1}(M;\mathbb{R}) \longrightarrow H^{1}(M^{1};\mathbb{R})$ is injective. By Lemma \ref{Lemma : Pullback}, the \'{A}lvarez class of $(M,\mathcal{F})$ vanishes. By Theorem \ref{Theorem : Alvarez Lopez} of \'{A}lvarez L\'{o}pez, $(M,\mathcal{F})$ is minimizable.
\end{proof}

\section{Rigidity of polycyclic groups}

We recall the definition of polycyclic groups.
\begin{defn}\label{Definition : Polycyclic groups}
A group $\Gamma$ is polycyclic if there exists a finite set $\{\Gamma_{i}\}_{i=0}^{n}$ of subgroups of $\Gamma$ such that
\begin{enumerate}
\item $\Gamma_{0}=\Gamma$, $\Gamma_{n}=\{1\}$,
\item $\Gamma_{i-1} \triangleright \Gamma_{i}$ and
\item $\Gamma_{i}/\Gamma_{i+1}$ is cyclic.
\end{enumerate}
\end{defn}
Clearly a polycyclic group is solvable. We refer to \cite{Segal} for general properties of polycyclic groups.

We show the following partial generalization of Lemma 2.1 in \cite{Meigniez}:
\begin{lem}\label{Lemma : Polycyclic groups1}
Let $G$ be a Lie group and $\Gamma$ be a dense polycyclic subgroup of $G$. We put $\Aut(G,\Gamma) = \{ a \in \Aut(G) \ | \ a(\Gamma)=\Gamma \}$. Let $\rho \colon \Aut(G,\Gamma) \longrightarrow \Aut(\Lie(G))$ be the action defined by $\rho(a)=Da$, where $Da$ is the restriction of the tangent map of $a \colon G \longrightarrow G$ to the Lie algebra of right invariant vector fields on $G$. Then every eigenvalue of $Da$ is an algebraic integer for every $a$ in $\Aut(G,\Gamma)$.
\end{lem}

\begin{proof}
We define $\Gamma_{0}=\Gamma$ and $\Gamma_{i}=[\Gamma_{i-1},\Gamma_{i-1}]$ for $i \geq 1$. Then there exists $n_{0}$ such that $\Gamma_{n_{0}}=\{1\}$ by the solvability of $\Gamma$. We define $G_{0}=G$ and $G_{i}=[G_{i-1},G_{i-1}]$ for $i \geq 1$. Since $\Gamma$ is dense in $G$, $\Gamma_{i}$ is dense in $G_{i}$ for every $i$. Then $G_{n_{0}}=\{1\}$ because $\Gamma_{n_{0}}=\{1\}$. We denote the canonical projection $G_{i} \longrightarrow G_{i}/G_{i+1}$ by $p_{i}$. Since every subgroup of a polycyclic group is polycyclic (see Exercise 2 of \cite{Segal}), $\Gamma_{i}$ is polycyclic for every $i$. In particular, $\Gamma_{i}$ is finitely generated. Then $p_{i}(\Gamma_{i})$ is finitely generated.

Let $a$ be an element of $\Aut(G,\Gamma)$. Then $a$ induces an automorphism $a_{i} \colon G_{i}/G_{i+1} \longrightarrow G_{i}/G_{i+1}$ for every $i$. Moreover $a_{i}$ preserves $p_{i}(\Gamma_{i})$. Note that the set of eigenvalues of $a$ is the union of the sets of eigenvalues of $a_{i}$. Then Lemma \ref{Lemma : Polycyclic groups1} follows from Lemma \ref{Lemma : Polycyclic groups 2} shown in below.
\end{proof}

\begin{lem}\label{Lemma : Polycyclic groups 2}
Let $\Gamma$ be a finitely generated dense subgroup of $\mathbb{R}^{d}$. Then every eigenvalue of any element of $\Aut(\mathbb{R}^{d},\Gamma)$ is an algebraic integer.
\end{lem}

\begin{proof}
Since $\Gamma$ is a finitely generated abelian group, $\Gamma \otimes_{\mathbb{Z}} \mathbb{R}$ is a finite dimensional Euclidean space. Let $a$ be an element of $\Aut(\mathbb{R}^{d},\Gamma)$. Then $a|_{\Gamma} \colon \Gamma \longrightarrow \Gamma$ is extended to $\tilde{a} \colon \Gamma \otimes_{\mathbb{Z}} \mathbb{R} \longrightarrow \Gamma \otimes_{\mathbb{Z}} \mathbb{R}$. We have a homomorphism $\pi \colon \Gamma \otimes_{\mathbb{Z}} \mathbb{R} \longrightarrow \mathbb{R}^{d}$ whose restriction to $\Gamma$ is equal to the injection $\Gamma \longrightarrow \mathbb{R}^{d}$. Since $\Gamma$ is dense in $\mathbb{R}^{d}$, $\pi$ is surjective. We have the following diagram:
\begin{equation}\label{Equation : Diagram1}
\xymatrix{ 
\Gamma \otimes_{\mathbb{Z}} \mathbb{R} \ar[dd]_{\pi} \ar[rrr]^{\tilde{a}} & & & \Gamma \otimes_{\mathbb{Z}} \mathbb{R}  \ar[dd]^{\pi} \\
         & \Gamma \ar[ul] \ar[ld] \ar[r]^{a} & \Gamma \ar[ur] \ar[rd] \\
\mathbb{R}^{d} \ar[rrr]_{a} & & & \mathbb{R}^{d}. }
\end{equation}
For a complex number $\lambda$, $\tilde{a} - \lambda \id_{\Gamma \otimes_{\mathbb{Z}} \mathbb{R}}$ is not surjective if $a - \lambda \id_{\mathbb{R}^{d}}$ is not surjective. Hence the set of the eigenvalues of $a$ is the subset of the eigenvalues of $\tilde{a}$. Since $\tilde{a}$ is regarded as an element of $\SL(n_{1};\mathbb{Z})$, where $n_1=\dim \Gamma \otimes_{\mathbb{Z}} \mathbb{R}$, the eigenvalues of $\tilde{a}$ are algebraic integers. The proof is completed.
\end{proof}

\section{Proof of Theorem \ref{Theorem : Rigidity}}

Let $\xi(\mathcal{F})$ and $\xi(\mathcal{F}^{1})$ be the \'{A}lvarez classes of $(M,\mathcal{F})$ and $(M^{1},\mathcal{F}^{1})$ respectively. Let $\gamma$ be a closed path in $M^{1}$. Then we have
\begin{equation}
\int_{\pi \circ \gamma} \xi(\mathcal{F}) = \int_{\gamma} \pi^{*} \xi(\mathcal{F}) = \int_{\gamma} \xi(\mathcal{F}^{1})
\end{equation}
by Lemma \ref{Lemma : Pullback}. Hence to show that $\int_{\pi \circ \gamma} \xi(\mathcal{F})$ is an algebraic integer, it suffices to show that $\int_{\gamma} \xi(\mathcal{F}^{1})$ is an algebraic integer. Note that every closed path in $M$ can be written as $\pi \circ \gamma$ for some closed path in $M^{1}$. Hence the proof of Theorem \ref{Theorem : Rigidity} is reduced to the case of transversely parallelizable foliations. 

By Proposition \ref{Proposition : Diagram}, the image of the holonomy homomorphism of the determinant bundle of the Molino's commuting sheaf of $(M^{1},\mathcal{F}^{1})$ is contained in the set of algebraic integers. Then, by Theorem \ref{Theorem : Alvarez Lopez 2} of \'{A}lvarez L\'{o}pez, the proof of Theorem \ref{Theorem : Rigidity} is completed.

\section{Proof of the Corollaries}

We deduce Corollary \ref{Corollary : Polycyclic} from Theorem \ref{Theorem : Rigidity}.
\begin{proof}[Proof of Corollary \ref{Corollary : Polycyclic}]
By Theorem \ref{Theorem : Rigidity}, it suffices to show that if $\pi_{1}M$ is polycyclic, then $\hol_{F}(\pi_{1}(F,x_{0}))$ is polycyclic. By the homotopy exact sequence of the fiber bundle $M^{1} \longrightarrow M$, $\pi_{1}M^{1}$ is an extension of $\pi_{1}M$ by an abelian group. Since the extension of a polycyclic group by a polycyclic group is polycyclic, $\pi_{1}M^{1}$ is polycyclic. By the homotopy exact sequence of the fiber bundle $M^{1} \longrightarrow W$, $\pi_{1}F$ is an extension of a subgroup $K$ of $\pi_{1}M^{1}$ by an abelian group. Since a subgroup of a polycyclic group is polycyclic (see Exercise 2 of \cite{Segal}), $K$ is polycyclic. Then $\pi_{1}F$ is polycyclic. Since the quotient of a polycyclic group by a normal subgroup is polycyclic, $\hol_{F}(\pi_{1}(F,x_{0}))$ is polycyclic. The proof is completed.
\end{proof}

We deduce Corollary \ref{Corollary : Nilpotent} from a theorem of Carri\`{e}re \cite{Carriere} and Theorem \ref{Theorem : Rigidity}. 
\begin{proof}[Proof of Corollary \ref{Corollary : Nilpotent}]
By Theorem \ref{Theorem : Rigidity}, it suffices to show that if $(M,\mathcal{F})$ is of polynomial growth, then $\hol_{F}(\pi_{1}(F,x_{0}))$ is polycyclic. By a theorem of Carri\`{e}re \cite{Carriere}, $(M,\mathcal{F})$ is of polynomial growth if and only if the structural Lie algebra of $(M,\mathcal{F})$ is nilpotent. Then $G$ is nilpotent by the assumption. Hence $\hol_{F}(\pi_{1}(F,x_{0}))$ is nilpotent. Since a finitely generated nilpotent group is polycyclic (see Corollary 8 of \cite{Segal}), $\hol_{F}(\pi_{1}(F,x_{0}))$ is polycyclic.
\end{proof}

\section{Examples}

\subsection{Riemannian foliations with semisimple structural Lie algebras}

Let $H_{1}$ be a semisimple Lie group and $H_{2}$ be a Lie group. Assume that we have a surjective homomorphism $\pi \colon H_{1} \longrightarrow H_{2}$. By Borel's lemma \cite{Borel}, we have a uniform lattice $L$ in $H_{1}$. Let $\mathcal{G}$ be the foliation defined by the fibers of $\pi$ on $H_{1}$. $\mathcal{G}$ induces a foliation $\mathcal{F}$ on $L \backslash H_{1}$. $\mathcal{F}$ is an $H_{2}$-Lie foliation. In particular, $(L \backslash H_{1},\mathcal{F})$ is transversely parallelizable, and therefore Riemannian. The structural Lie algebra of $(L \backslash H_{1},\mathcal{F})$ is the Lie algebra of the closure $K$ of $\pi(L)$ in $H_{2}$. Since the quotient of a semisimple Lie group is semisimple, $H_{2}$ is semisimple. Then $K$ is semisimple. By Theorem \ref{Theorem : Reductive}, $(L \backslash H_{1},\mathcal{F})$ is minimizable.

\subsection{Carri\`{e}re's example}

We present Carri\`{e}re's example which is a Riemannian foliation constructed in a similar way to the previous one. Let $A$ be a matrix in $\SL(2;\mathbb{Z})$ with two positive real eigenvalues. We can define $A^{t}$ for real numbers $t$.

Let $\mathbb{R} \ltimes_{A} \mathbb{R}^{2}$ be the Lie group with the underlying manifold $\mathbb{R} \times \mathbb{R}^{2}$ and the product defined by $(t,x) \cdot (t',x') = (t + t', A^{t'} x + x')$. Let $\mathbb{Z} \ltimes_{A} \mathbb{Z}^{2}$ be the subgroup of $\mathbb{R} \ltimes_{A} \mathbb{R}^{2}$ with the underlying set $\mathbb{Z} \times \mathbb{Z}^{2}$. Let $v_1$ and $v_2$ be eigenvectors of $A$. We denote the eigenvalues of $A$ with respect to $v_i$ by $\lambda_i$ for $i=1$ and $2$. Let $\mathbb{R} \ltimes_{\lambda_1} \mathbb{R}$ be the Lie group with the underlying manifold $\mathbb{R} \times \mathbb{R}$ and the product defined by $(t,u) \cdot (t',u') = (t + t', \lambda_{1}^{t'} u + u')$. We define a map
\begin{equation}
\begin{array}{cccc}
\pi \colon & \mathbb{R} \ltimes_{A} \mathbb{R}^{2} & \longrightarrow & \mathbb{R} \ltimes_{\lambda_1} \mathbb{R} \\
 & (t, u_{1} v_{1}+ u_{2} v_{2}) & \longmapsto & (t, u_{1} v_{1})
\end{array}
\end{equation}
for $u_{1}$, $u_{2} \in \mathbb{R}$. Then $\pi$ is a surjective homomorphism. The fibers of $\pi$ induces a $1$-dimensional foliation $\mathcal{F}$ on $T_{A}^{3}=(\mathbb{Z} \ltimes_{A} \mathbb{Z}^{2}) \backslash (\mathbb{R} \ltimes_{A} \mathbb{R}^{2})$. $(T_{A}^{3},\mathcal{F})$ is an $(\mathbb{R} \ltimes_{\lambda_1} \mathbb{R})$-Lie foliation. Carri\`{e}re showed that $(T_{A}^{3},\mathcal{F})$ is not minimizable by computing its basic cohomology (see \cite{Carriere 2}).

$T_{A}^{3}$ is a $T^2$-bundle over $S^1$ and the leaves of $\mathcal{F}$ are contained in the $T^2$-fibers. The restriction of $\mathcal{F}$ to each $T^2$-fiber is a linear flow with dense leaves. We will compute the \'{A}lvarez class of $(T_{A}^{3},\mathcal{F})$ in the next example.

\subsection{Torus bundles over $S^1$}

Let $A$ be an element of $\SL(n;\mathbb{Z})$ with an eigenvector $v$ with respect to an eigenvalue $\lambda$. Assume that $\lambda$ is a positive real number and the components of $v$ are real numbers linearly independent over $\mathbb{Q}$. $A$ induces a diffeomorphism $\overline{A}$ of $T^{n}=\mathbb{R}^{n}/\mathbb{Z}^{n}$. We denote the mapping torus $[0,1] \times T^{n}/(0,\overline{A}w) \sim (1,w)$ of $\overline{A}$ by $M$. We define a map $\pi \colon M \longrightarrow S^1$ by $\pi([(t,w)])=t$. Then $\pi$ defines a $T^{n}$-bundle over $S^1$. Since $v$ is preserved by $A$ up to multiplication of a real number, we have a $1$-dimensional foliation $\mathcal{F}$ on $M$ defined by the lines parallel to $v$ in each $T^{n}$-fiber of $\pi$. By the assumption on $v$, the leaves of $\mathcal{F}$ are dense in the fibers of $\pi$. Since we can define a bundle-like metric $g$ of $(M,\mathcal{F})$ by $g = \rho(t)g_0 + (1-\rho(t))\overline{A}^{*}g_0 + dt \otimes dt$, where $g_0$ is the flat metric on $T^n$ and $\rho(t)$ is a smooth non-negative function which satisfies $\rho(t)=0$ near $0$ and $\rho(t)=1$ near $1$, $(M,\mathcal{F})$ is Riemannian. Note that the structure Lie algebra of $(M,\mathcal{F})$ is abelian.

Since the linear flows on tori are minimizable, the restriction of the \'{A}lvarez class on the fiber of $\pi$ is trivial. Let $\gamma$ be a closed path in $M$ such that $\pi \circ \gamma$ is a representative of the generator of the fundamental group of $S^1$ and the holonomy of $\pi$ with respect to $\pi \circ \gamma$ is $A$. The holonomy of the Molino's commuting sheaf of $(M,\mathcal{F})$ is induced by $A$. By Proposition \ref{Proposition : Diagram} and Theorem \ref{Theorem : Alvarez Lopez 2} of \'{A}lvarez L\'{o}pez. we have
\begin{equation}
\int_{\gamma}\xi(\mathcal{F})= \log \frac{1}{\lambda},
\end{equation}
where $\xi(\mathcal{F})$ is the \'{A}lvarez class of $(M,\mathcal{F})$.

\subsection{A Riemannian foliation on a solvmanifold}

We present an example of a $2$-dimensional Riemannian foliation on a $6$-dimensional solvmanifold with nonabelian nilpotent structure Lie algebra and nontrivial \'{A}lvarez class.

Let $p$ be a prime and $\alpha$ be an element of $\mathbb{Z}(\sqrt{p})$ which has an inverse $\beta$ in $\mathbb{Z}(\sqrt{p})$. Set $k=\GCD(\alpha_2,\beta_2)$, where $\alpha_2$ and $\beta_2$ are integers satisfying $\alpha=\alpha_1+\alpha_2\sqrt{p}$ and $\beta=\beta_1+\beta_2\sqrt{p}$ for some integers $\alpha_1$ and $\beta_1$. Let $G$ be the nilpotent Lie group defined by
\begin{equation}
G=\Big\{ \begin{pmatrix} 1 & x & z \\ 0 & 1 & y \\ 0 & 0 & 1 \end{pmatrix} \ \Big| \ x,y,z \in \mathbb{R} \Big\}
\end{equation}
and $\Gamma$ the subgroup of $G$ generated by
\begin{equation}
\begin{array}{c}
A_1=\begin{pmatrix} 1 & 1 & 0 \\ 0 & 1 & 0 \\ 0 & 0 & 1 \end{pmatrix}, \quad
A_2=\begin{pmatrix} 1 & 0 & 0 \\ 0 & 1 & 1 \\ 0 & 0 & 1 \end{pmatrix}, \\ \vspace{-8pt} \\
A_3=\begin{pmatrix} 1 & k\sqrt{p} & 0 \\ 0 & 1 & 0 \\ 0 & 0 & 1 \end{pmatrix}, \quad
A_4=\begin{pmatrix} 1 & 0 & 0 \\ 0 & 1 & k\sqrt{p} \\ 0 & 0 & 1 \end{pmatrix}.
\end{array}
\end{equation}
$\Gamma$ is dense in $G$.

Consider the Lie group $H=G \oplus \mathbb{R}^{2}$. Let $\iota'$ be a homomorphism $\iota' \colon \Gamma \longrightarrow \mathbb{R}^{2}$ which is defined by
\begin{equation}
\iota'\Big(\begin{pmatrix} 1 & x_1 + x_2 \sqrt{p} & z \\ 0 & 1 & y_1 + y_2 \sqrt{p} \\ 0 & 0 & 1 \end{pmatrix} \Big)=(x_2,y_2),
\end{equation}
where $x_1,x_2,y_1$ and $y_2$ are integers, and we define an embedding $\iota \colon \Gamma \longrightarrow H$ by $\iota(g)=(g,\iota'(g))$ for every $g$. Then $\iota(\Gamma)$ is a uniform lattice of $H$. The fibers of the first projection $H \longrightarrow G$ are preserved by the right multiplication of elements of $\iota(\Gamma)$ and define a $G$-Lie foliation $\mathcal{G}$ of dimension $2$ and codimension $3$ on $H/\iota(\Gamma)$.

Let $\begin{pmatrix} a' & b' \\ c' & e' \end{pmatrix}$ be an element of $\SL(2;\mathbb{Z})$ and set $\begin{pmatrix} a & b \\ c & e \end{pmatrix}=\begin{pmatrix} a'\alpha & b'\alpha \\ c'\alpha & e'\alpha \end{pmatrix}$. Let $f$ be the map $G \longrightarrow G$ defined by 
\begin{equation}
f\begin{pmatrix} 1 & x & z \\ 0 & 1 & y \\ 0 & 0 & 1 \end{pmatrix}=\begin{pmatrix} 1 & ax+by & \alpha^{2} z+ \frac{1}{2}acx^2 + bcxy + \frac{1}{2}bey^2 \\ 0 & 1 & cx+ey \\ 0 & 0 & 1 \end{pmatrix}.
\end{equation}
Then $f$ is an automorphism of $G$ whose restriction to $\Gamma$ is an automorphism of $\Gamma$ by the definition of $k$ and $\Gamma$.

Since $\iota(\Gamma)$ is a uniform lattice of $H$, the automorphism $f|_{\Gamma}$ of $\Gamma$ is uniquely extended to an automorphism $\tilde{f}$ of $H$ and induces a diffeomorphism $\overline{f}$ of $H/\iota(\Gamma)$ by Mal'cev theory \cite{Malcev} (see also \cite{Raghunathan}). Since $\overline{f}$ preserves $\mathcal{G}$, we have a foliation $\mathcal{F}$ on the mapping torus $M$ of $\overline{f}$, which is Riemannian. We have $f^{*}(dx \wedge dy \wedge dz)=\alpha^{3} dx \wedge dy \wedge dz$. By Proposition \ref{Proposition : Diagram} and Theorem \ref{Theorem : Alvarez Lopez 2} of \'{A}lvarez L\'{o}pez, we have 
\begin{equation}
\int_{S^1}\xi(\mathcal{F})= \log \alpha^3
\end{equation}
where $S^1$ is the base space of the canonical fibration $M \longrightarrow S^1$ of the mapping torus and $\xi(\mathcal{F})$ is the \'{A}lvarez class of $(M,\mathcal{F})$.

\subsection{Meigniez's examples}

Let $G$ be a solvable Lie group. Let $\Gamma$ be a finitely generated subgroup $\Gamma$ of $G$ which contains a uniform subgroup of $G$. Meigniez \cite{Meigniez} showed that there exists a $G$-Lie foliation $\mathcal{F}$ on a close manifold such that the image of the holonomy homomorphism of $(M,\mathcal{F})$ is $\Gamma$. By Theorem \ref{Theorem : Alvarez Lopez 2} of \'{A}lvarez L\'{o}pez, the period map of the \'{A}lvarez class of a $G$-Lie foliation on a closed manifold $M$ coincides with the composite
\begin{equation}
\xymatrix{ \pi_{1}M \ar[r]^>>>>>{\hol_{M}} & G \ar[r]^>>>>>{\modular} & \mathbb{R}_{>0} \ar[r]^>>>>>{\log} & \mathbb{R} }
\end{equation}
where $\hol_{M}$ is the holonomy homomorphism of $(M,\mathcal{F})$ and $\modular$ is the modular function of the Lie group $G$. See \cite{Meigniez} for a one-parameter family of $GA(1)$-Lie foliations on a closed manifold whose fundamental group is not polycyclic such that the periods of the \'{A}lvarez classes change nontrivially with respect to the parameter.

\end{document}